\documentclass{article}
\usepackage{graphicx} 
\graphicspath{ {./images/} }

\usepackage{float}
\usepackage{amsmath}
\usepackage{amssymb}
\usepackage{graphicx,psfrag}
\usepackage{amsfonts,amsthm}
\usepackage{amsxtra,enumerate}

\def\veps{\varepsilon}

\theoremstyle{plain}
\newtheorem{theorem}{Theorem}[section]

\theoremstyle{definition}
\newtheorem{definition}[theorem]{Definition}
\theoremstyle{remark}

\newtheorem{example}[theorem]{Example}

\newcommand{\pd}[2]{\frac{\partial^2{#1}}{\partial{#2}^2 }}

\newcommand{\R}{\ensuremath{\mathbb{R}}}
\newcommand{\N}{\ensuremath{\mathbb{N}}}
\newcommand{\Z}{\ensuremath{\mathbb{Z}}}
\newcommand{\T}{\ensuremath{\mathbb{T}}}
\newcommand{\C}{\ensuremath{\mathbb{C}}}
\newcommand{\deld}{^{\Delta}}
\newcommand{\crd}{{\rm C}_{\rm rd}}
\newcommand{\dd}{{\rm d}}

\title{Solving the Wave Equation on Discrete Time Scales}
\author{Davis Funk and Charis Tsikkou}
\date{\today}

\begin{document}

\maketitle

\begin{abstract}
\noindent This paper presents a solution to an initial value problem for the 1-dimensional wave equation on time scales through the application of a Fourier transform and its inverse via contour integrals. The time scale of the spatial dimension is set to the integers and a broader class of discrete time scales, while the time dimension is set to the positive real numbers.
\end{abstract}


\noindent


\section{Introduction}
\quad Time scales calculus is a relatively new field of mathematics introduced in 1988 by Stefan Hilger \cite{hilger}. This field unifies concepts from the traditional continuous calculus and discrete mathematics, and creates generalized definitions that account for both cases - as well many arbitrary cases in between. Many important mathematical objects, such as the derivative and integral, have well-studied time scales definitions. When the time scale is chosen to be the real numbers or the integers, then these definitions reduce to the equivalent objects from calculus and discrete math. These time scales definitions can be applied in fields like differential equations, which is known as dynamic equations when dealing with derivatives on time scales. A variety of other fields, such as control theory, special functions, and fractional calculus have been studied in depth in a time scales framework \cite{besselpaper,controlpaper,fractionalpaper}. 

\quad Partial differential equations, and multi-dimensional objects in general, are still under development \cite{multivariate}. In recent work, a time scales univariate Laplace transform was used to solve the wave equation on time scales \cite{jacksonpaper}. When the wave equation was previously studied, a general convolution theorem was not known at the time, leading to a limited class of functions being considered for initial data, limiting the number of possible solutions. The goal of this paper is to solve the wave equation with the spatial variable set to an arbitrary isolated time scale, and the time variable lying in $(0,\infty)$. Cuchta and Ferriera recently developed a novel formulation of the Fourier transform that utilizes a contour integral \cite{heatpaper}. This formulation was then used to solve the heat equation, and is applied by us in this paper in order to solve the wave equation.

On time scales, we define the wave equation based on Jackson's definition \cite{jacksonpaper}:
\begin{equation}\label{Eq1.1}
    u^{\Delta^2_{tt}}\ = \ c^2 u^{\Delta^2_{xx}}.
\end{equation}

In this paper, we first solve the following wave equation problem:
\begin{equation}\label{Eq1.2}
\begin{cases}
\pd{u}{t} = \ c^2 u^{\Delta^2_{xx}}, t\in [ 0,\infty), x\in\T = \Z \\
u(0,x) = f(x) \\
u_t(0,x) = g(x)
\end{cases}
\end{equation}

\section{Preliminaries and Definitions}
\quad This section offers a short introduction to the relevant concepts from time scales calculus, as well as to the formulations for the Fourier transform and its inverse that will be used later.

\begin{definition}
A time scale, $\T$, is an arbitrary, nonempty closed subset of the real numbers.
\end{definition} 
Some common examples of time scales include:
\begin{itemize}
    \item The real numbers - $\R$
    \item The integers - $\Z$
    \item The quantum numbers - $q^{\N_0} = \lbrace q^{n} | n \in \N_0 \rbrace$, where $q>1$.
\end{itemize}
Numerous important operators are defined in time scales calculus, upon which definitions of familiar mathematical objects are built.

\begin{definition} The following three operators are key for further developing time scales calculus.
\begin{itemize}
\item
The forward jump operator $\sigma\colon\T\to \T$ is defined as $\sigma(t) := \inf\lbrace s\in\T \ | \ s>t\rbrace$.
\item
The backward jump operator $\rho\colon\T\to\T$ is defined as $\rho(t) := \sup\lbrace s\in\T \ | \ s<t\rbrace$.
\item
The graininess operator - $\mu\colon\T\to [0,\infty)$ is defined by $\mu(t) := \sigma(t)-t$.
\end{itemize}
\end{definition}

\begin{example}
When we work on the real numbers, $\sigma(t)$ is simply $t$ itself. Due to this, $\mu$ is always zero on the time scale of real numbers.
When we work with integers, $\sigma(t)$ is $t+1$, and $\mu$ is identically equal to one at every point. 
\end{example}
Next, we use these operators to define the time scales derivative. When working with time scales derivatives we frequently use the notation $\T^\kappa$, where $\T^{\kappa}=\T\setminus\{\max \T\}$ given $\max \T$ exists; otherwise $\T^{\kappa}=\T$.
\begin{definition}
Suppose $f\colon\T\to\R$ and let $t\in\T^\kappa$. We say $f$ is delta differentiable at $t\in\T^\kappa$ if there exists an $\alpha$ such that for any $\veps>0$, there is a neighborhood $U$ of $t$ such that
\begin{equation*}
\left|[f(\sigma(t))-f(s)]-\alpha[\sigma(t)-s]\right| \leq
\veps|\sigma(t)-s|\quad \text{for all }  \quad s \in U.
\end{equation*}
\end{definition}
We denote this $\alpha$ by $f\deld(t)$ and call it the delta derivative of $f$ at t. The following theorem from Agarwal et al. \cite{tsbasics} can be used to easily calculating the delta derivative, depending on how the time scale is defined.

\begin{theorem}
If $f$ is delta differentiable at $t\in\T^\kappa$, then it follows that
\begin{itemize}
\item[a.] If $\mu(t) > 0$, then $f\deld(t)=\frac{f(\sigma(t))-f(t)}{\mu(t)}.$
\item[b.] If $\mu(t) = 0$, then $f\deld(t)=\lim_{s \to t}\frac{f(t)-f(s)}{t-s}.$
\end{itemize}
\end{theorem}
Next, we talk about important concepts regarding integrals on time scales.
\begin{definition}
A function $f\colon\T\to\R$, is considered rd-continuous if it is continuous at every right dense point and if the left sided limit exists at every left dense point.
\end{definition}
A right dense point is a point with at which $\sigma(t) = t$, and a left dense point is one where $\rho(t) = t$. A point that is neither left or right dense is called an isolated point. The set of all rd continuous functions $f\colon\T\to\R$ is denoted by $\crd$. The following theorem from Agarwal et al. \cite{tsbasics} shows why this property is valuable.
\begin{theorem}
Every rd-continuous function, $f\colon\T\to\R$, has an antiderivative $F\colon\R\to\T$, such that $F\deld(t) = f$ on all $t\in\T^{\kappa}$
\end{theorem}
Now, we define the time scale integral as follows:
\begin{definition}
Let f be an rd-continuous function. Then for all $a,b \in \T$, the Cauchy integral of f is defined as:
\begin{equation*}
\int^{b}_{a}f(t)\Delta t=F(b)-F(a), \ \ \ \forall a,b\in\T.
\end{equation*}

\end{definition}
Depending on the time scale, this integral takes different forms. The following theorem from Bohner and Peterson \cite{greenbook} shows us how the time scales integral changes:
\begin{theorem}\label{ex2}
Let $a,b\in\T$ with $a<b$ and assume that $f\in\crd$.
\begin{itemize}
\item[a.] When $\T=\R$, then
$$ \int^{b}_{a}f(t)\Delta t=\int^{b}_{a}f(t)\dd t. $$
\item[b.] If the interval $[a,b]\cap\T$ only consists of isolated points
$$ \int^{b}_{a}f(t)\Delta t=\sum_{t\in[a,b)}\mu(t)f(t). $$
\item[c.] When $\T=h\Z=\lbrace hk: \ k\in\Z\rbrace$ for $h>0$, then
$$ \int^{b}_{a}f(t)\Delta t=h\sum^{b/h-1}_{k=a/h}f(kh). $$
\end{itemize}
\end{theorem}
The Fourier transform is helpful when solving a time scales partial dynamic equation. An important property for a function to have, in order for the Fourier transform to be defined is for that function to be regulated.
\begin{definition}
    A function $f\colon\T\to\R$ is regulated if its finite, right-sided limits exist for all right dense points in \T, and if its finite, left-side limits exist for all left dense points in \T.
\end{definition}
For our solution we use the Fourier transform definition given by Cuchta and Georgiev \cite{bilaterallap}. Before defining this, we must discuss the time scales exponential and the cylinder transform.
\begin{definition}
$\xi_{h}(z)$ is the cylinder transform, defined as $$\xi_{h}(z) =
\begin{cases} 
      \frac{Log(1+zh)}{z} & h\neq 0 \\
      z & h = 0 \\
\end{cases}$$
\end{definition}
We define the time scales exponential as
\begin{definition}
\begin{equation*}
    e_{p}(t,s) \ = \exp \biggl\{ \int^{t}_{s}\xi_{\mu(\tau)}(p(\tau))\Delta\tau \biggr\}.
\end{equation*}
\end{definition}
For this definition, p must be rd-continuous and regressive. \\
\begin{definition}
A function $p\colon\T\to\R$ is regressive if $1+\mu(t)p(t)\neq0$ for all $t\in\T$
\end{definition}
Also important to note is the circle minus operator $\ominus$.
For a function $p$, we define $\ominus p$ as
$\ominus p \ = \ \frac{-p}{1+\mu p} \\$
\\
Now, we define the transform as follows:
\begin{definition}
Suppose $f\colon\T\to\R$ is regulated. Then the Fourier transform centered at $s\in\T$ is defined as:
\begin{equation}\label{ft}
    \mathcal{F}_{\T}\lbrace f \rbrace (z;s) \ = \int^{\infty}_{-\infty}f(t)e_{\ominus iz}(\sigma(t),s)\Delta t,
\end{equation}
whenever the integral converges.
\end{definition}
Cuchta and Georgiev \cite{bilaterallap} also prove the following theorem:
\begin{theorem}
If $f\colon\T\in\C$ is n-times delta differentiable and for all $n \ \epsilon \ \lbrace0,....,k-1\rbrace, \lim_{t \to \pm\infty}f^{\Delta^{k}}(t)e_{\ominus iz}(t,s)=0$, then
\begin{equation}
    \mathcal{F}_{\T} \biggl\{ f^{\Delta^{k}} \biggr\} (z;s) \ = \ (iz)^k \mathcal{F}_{\T}\biggl\{ f \biggr\} (z;s)
\end{equation}
for all $z\in\C$ in which the respective integrals exist.
\end{theorem} 
When we set $\T=\Z$ and $s=0$ we find that equation \ref{ft} simplifies to the form
\begin{equation}
    \mathcal{F}_{\T} \lbrace f \rbrace (z;0) \ = \ \sum_{k\in\Z} \frac{f(k)}{(1+iz)^{k+1}}.
\end{equation}
We can invert this version of the Fourier transform with the following theorem defining an inverse Fourier transform as shown by Cuchta and Ferreira \cite{heatpaper}.
\begin{theorem}
If \ $\T=\Z$ and F is analytic, except possibly on a branch cut $(-\infty,0]$, then
\begin{equation*}
    \mathcal{F}^{-1} \lbrace F \rbrace (t;0) \ = \ \frac{1}{2\pi} \oint_{C} F(z)(1+iz)^t \mathrm{d} z,
\end{equation*}
where C is a circle with center $z = i$ and radius $0 < r < 1$.
\end{theorem}

We can also look at a broader class of discrete time scales, where our graininess function is injective. In this situation, when solving the contour integral introduced by the inverse Fourier transform, we only need to deal with simple poles \cite{complexbook}, whose residues are straightforward to calculate. A similar procedure can be done for time scales with a non injective graininess as well.

To define the inverse Fourier transform on this class of time scales, we state the following theorem from Cuchta and Ferreira \cite{heatpaper}.

\begin{theorem}
Let $\T$ be a countable time scale of the form $\lbrace...,t_{-1},t_0,t_1,...\rbrace$, where for all integers a and b, $t_a < t_b$ if and only if $a < b$. If F is analytic, except possibly on the branch cut $(-\infty,0]$, and the graininess function $\mu: \T \in [0,\infty)$ associated with the time scale is injective, then 
\begin{equation*}
    \mathcal{F}^{-1} \lbrace F \rbrace (t;s) \ = \ \frac{1}{2\pi} \oint_{C} F(z)e_{iz}(t,s) \dd z
\end{equation*}
where $\oint_{C}$ is defined as $\lim_{m \to \infty}\oint_{C_m}$, where $C_m$ is the contour containing all points $\biggl\{ \frac{i}{t_k}: -m \leq k \leq m \biggr\}$
\end{theorem}

Finally, we want to introduce the definition of convolution on time scales. To develop this definition we must discuss the shifting problem. We use Cuchta and Ferriera's definition in \cite{heatpaper}.
\begin{definition}
The following problem is known as the shifting problem: \\
\begin{equation}\label{shift}
\begin{cases}
u^{\Delta_t}(t,\sigma(s)) = -u^{\Delta_s}(t,s),  t,s \in \T, \\
u(t,t_0) = f(t), \quad t\in\T
\end{cases}
\end{equation}
where $f\colon[t_0,\infty)\cap\T\to\R$. The solution to this problem is denoted by $\hat{f}(t,s) = u(t,s)$
\end{definition}
We can do define the time scales convolution using the shift as developed in \cite{convolution}.
\begin{definition}
For given functions $f,g\colon\T\to\R$, their convolution, $f \ast g$ is defined by
\begin{equation*}
(f \ast g)(t) = \int_{t0}^t \hat{f}(t,\sigma(s)) g(s) \Delta s, \quad t\in\T,
\end{equation*}
where $\hat{f}$ is the shift discussed in \ref{shift}.
\end{definition}
\section{Solving the wave equation with $x\in\Z$}
On time scales, we define the wave equation based on equation \ref{Eq1.1}.
First, we solve the wave equation with the following conditions: $t\in[0,\infty), x\in\Z$. Under these conditions this equation becomes:
\begin{equation*}
    \pd{u}{t} = \ c^2 u^{\Delta^2_{xx}},
\end{equation*}
where we define initial conditions $u(0,x) = f(x)$ and $u_t(0,x) = g(x)$. \\
Like when solving a normal partial differential equation, we apply a Fourier transform to turn our partial dynamic equation into an ordinary dynamic equation \cite{pdebook}. We will also apply Theorem 3, which results in the following

\begin{equation*}
    \pd{\hat{u}}{t}(t,z) = \ c^2 (iz)^2 \hat{u}.
\end{equation*}
This transforms our initial conditions to $\hat{u}(0,z) = \hat{f}(z)$ and $\hat{u}_t(0,z) = \hat{g}(z)$.
We solve this second order ODE by looking for a solution of the form $\hat{u}(t,z) = e^{\lambda t}$. After plugging this into our ODE, and solving the resulting characteristic equation we find that $\lambda = \pm czi$. \\
Thus, our general equation for $\hat{u}$ is defined as:
\begin{equation*}
    \hat{u}(t,z) = c_1(z)\cos{(czt)} + c_2(z)\sin{(czt)}.
\end{equation*}
By inputting the initial condition $\hat{u}(0,z) = \hat{f}(z)$, we find that the first coefficient $c_1(z) = \hat{f}(z)$ To find $c_2(z)$, we calculate
\begin{equation*}
    \hat{u}_t(t,z) = -czc_1(z)\sin{(czt)} + czc_2(z)\cos{(czt)}.
\end{equation*}
By using the initial condition $\hat{u}_t(0,z) = \hat{g}(z)$ we find the second coefficient $c_2(z) = \frac{\hat{g}(z)}{cz}$.
Thus, our final expression for $\hat{u}$ is
\begin{equation*}
    \hat{u}(t,z) = \hat{f}(z)\cos{(czt)} + \frac{\hat{g}(z)}{cz}\sin{(czt)}.
\end{equation*} \\
To find $u(t,z)$, we compute the inverse Fourier transforms of the cosine and sine terms. First, define $G_{c,t}(z) = \cos{(czt)}$, and take the inverse transform

\begin{equation*}
    \mathcal{F}^{-1} \lbrace G_{c,t}(z) \rbrace (x;0) \ = \ \frac{1}{2\pi} \oint_{C} \cos{(czt)} \ (1+iz)^x \dd z
\end{equation*}
When $x\geq0$, the function in the integral is analytic, and our integral vanishes according to the Cauchy Goursat theorem \cite{complexbook}.
When $x<0$,
\begin{align*}
    \mathcal{F}^{-1} \lbrace G_{c,t}(z) \rbrace (x;0) \ & = \ \frac{1}{i^{|x|}}\frac{1}{2\pi} \oint_{C} \frac{\cos{(czt)}}{(z-i)^{|x|}} \mathrm{d} z \\ 
    & = \ \frac{1}{2\pi i^{|x|}} \ \frac{2\pi i}{(|x|-1)!} \frac{\mathrm{d}^{|x|-1}}{\mathrm{d}z^{|x|-1}} [\cos{(czt)}]_{z=i} \\
    & = \ \frac{1}{i^{|x|-1}} \ \frac{1}{(|x|-1)!} \frac{d^{|x|-1}}{dz^{|x|-1}} [\cos{(czt)}]_{z=i}.
\end{align*}
By repeatedly differentiating the cosine, we see:
\begin{itemize}
\item If $|x|$ odd, with $|x|=2k+1$, then
\begin{equation*}
\frac{\mathrm{d}^{2k}}{\mathrm{d}z^{2k}}[\cos{(czt)}]_{z=i} = (-1)^k c^{2k} t^{2k} \cos{(cit)}. \
\end{equation*}
\item If $|x|$ even, with $|x|=2k$, then
\begin{equation*}
\frac{\mathrm{d}^{2k-1}}{\mathrm{d}z^{2k-1}}[\cos{(czt)}]_{z=i} = (-1)^{k} c^{2k-1} t^{2k-1} \sin{(cit)}.
\end{equation*}
\end{itemize}

To simplify these expressions, we can use the following identities:
\begin{equation*}
    \sin(cit) \ = \ \sin{(0)}\cosh{(ct)}+i\cos{(0)}\sinh{(ct)} \ = \ i\sinh{(ct)},
\end{equation*}
\begin{equation*}
    \cos(cit)  \ = \ \cos{(0)}\cosh{(ct)}-i\sin{(0)}\sinh{(ct)} \ = \ \cosh{(ct)}.
\end{equation*}
Thus we write,
\begin{equation*}
\mathcal{F}^{-1} \lbrace G_{c,t}(z) \rbrace (x;0) \ =
\begin{cases} 
      \frac{1}{i^{2k+1-1}}\frac{1}{(2k+1-1)!}(-1)^{k} c^{2k} t^{2k} \cosh{(ct)}, \quad & |x| = 2k+1 \\ \\
      \frac{1}{i^{2k-1}}\frac{1}{(2k-1)!}(-1)^{k} c^{2k-1} t^{2k-1} i\sinh{(ct)}, \quad & |x| = 2k \\
\end{cases}
\end{equation*}
\begin{equation*}
\mathcal{F}^{-1} \lbrace G_{c,t}(z) \rbrace (x;0) \ =
\begin{cases} 
      \frac{(-1)^k c^{2k} t^{2k}}{i^{2k}(2k)!} \cosh{(ct)}, \quad & |x| = 2k+1 \\ \\
      \frac{(-1)^{k} c^{2k-1} t^{2k-1}}{i^{2k-1} (2k-1)!} i\sinh{(ct)}, \quad  & |x| = 2k \\
\end{cases}
\end{equation*}
Finally, we obtain,
\begin{equation*}
\mathcal{F}^{-1} \lbrace G_{c,t}(z) \rbrace (x;0) \ =
\begin{cases} 
      \frac{(ct)^{2k}}{(2k)!} \cosh{(ct)}, \quad & |x| = 2k+1 \\ \\
      \frac{-(ct)^{2k-1}}{(2k-1)!} \sinh{(ct)}, \quad & |x| = 2k \\
\end{cases}
\end{equation*}
Next, we calculate the inverse Fourier transform of the sine term. First, we define $H_{c,t} = \frac{\sin{(czt)}}{z}$, and we calculate
\begin{equation*}
    \mathcal{F}^{-1} \lbrace H_{c,t}(z) \rbrace (x;0) \ = \ \frac{1}{2\pi} \oint_{C} \frac{\sin{czt}}{z} \ (1+iz)^x \dd z
\end{equation*}
When $x\geq0$, the function in the integral is analytic, and our integral vanishes. When $x<0$,
\begin{align*}
    \mathcal{F}^{-1} \lbrace H_{c,t}(z) \rbrace (x;0) \ & = \ \frac{1}{2\pi i^{|x|}} \ \frac{2\pi i}{(|x|-1)!} \frac{\mathrm{d}^{|x|-1}}{\mathrm{d}z^{|x|-1}} \Bigg[\frac{\sin{(czt)}}{z}\Bigg]_{z=i} \\
    & = \ \frac{1}{i^{|x|-1}} \ \frac{1}{(|x|-1)!} \frac{\mathrm{d}^{|x|-1}}{\mathrm{d}z^{|x|-1}} \Bigg[\frac{\sin{(czt)}}{z}\Bigg]_{z=i}. \
\end{align*}
By repeatedly differentiating the sine term, we can find \\
For $|x|$ odd, with $|x|=2k+1$
\begin{equation*}
\frac{\mathrm{d}^{2k}}{\mathrm{d}z^{2k}}\Bigg(\frac{\sin{(czt)}}{z}\Bigg)_{z=i} = \ (2k)! \Bigg[ \sum^{k}_{m=0}\frac{(-1)^m (ct)^{2m}}{(2m)!(i)^{2k-2m+1}} \sin{(cit)}  -  \sum^{k-1}_{m=0}\frac{(-1)^m (ct)^{2m+1}}{(2m+1)!(i)^{2k-2m}} \cos{(cit)} \Bigg]
\end{equation*}
For $|x|$ even, with $|x|=2k$
{\small
\begin{equation*}
\frac{\mathrm{d}^{2k+1}}{\mathrm{d}z^{2k+1}}\Bigg(\frac{\sin{(czt)}}{z}\Bigg)_{z=i} = \ (2k+1)! \Bigg[ \sum^{k}_{m=0}\frac{-(-1)^m (ct)^{2m}}{(2m)!(i)^{2k-2m+2}} \sin{(cit)} + \sum^{k}_{m=0}\frac{(-1)^m (ct)^{2m+1}}{(2m+1)!(i)^{2k-2m+1}} \cos{(cit)} \Bigg] 
\end{equation*}}
Next, for $x<0$ we compute, \\
{\footnotesize
$$\mathcal{F}^{-1} \lbrace H_{c,t}(z) \rbrace(x;0) =
\begin{cases}
& \frac{(2k)!}{i^{2k+1-1}(2k+1-1)!} \Bigg[ \sum^{k}_{m=0}\frac{(-1)^m (ct)^{2m}}{(2m)!(i)^{2k-2m+1}} i\sinh{(ct)} \\
&-  \sum^{k-1}_{m=0}\frac{(-1)^m (ct)^{2m+1}}{(2m+1)!(i)^{2k-2m}} \cosh{(ct)} \Bigg], \ \ \ \ \ \text{when} \ \ |x| \  \ odd \\ \\
& \frac{(2k-1)!}{i^{2k-1}(2k-1)!} \Bigg[ \sum^{k-1}_{m=0}\frac{-(-1)^m (ct)^{2m}}{(2m)!(i)^{2k-2m}} i\sinh{(ct)}\\
&+ \sum^{k-1}_{m=0}\frac{(-1)^m (ct)^{2m+1}}{(2m+1)!(i)^{2k-2m-1}} \cosh{(ct)} \Bigg], \ \  \ \  \text{when} \ \ |x| \ even 
\end{cases}$$}
We simplify to reach \\
{\footnotesize
$$\mathcal{F}^{-1} \lbrace H_{c,t}(z) \rbrace (x;0) =
\begin{cases}
& \frac{1}{i^{2k}} \Bigg[ \sum^{k}_{m=0}\frac{(-1)^m (ct)^{2m}}{(2m)!(i)^{2k-2m+1}} i\sinh{(ct)}  \\ 
&-  \sum^{k-1}_{m=0}\frac{(-1)^m (ct)^{2m+1}}{(2m+1)!(i)^{2k-2m}} \cosh{(ct)} \Bigg], \ \ \ \text{when} \  \  |x| \ odd \\ 
 & \frac{1}{i^{2k-1}} \Bigg[ \sum^{k-1}_{m=0}\frac{-(-1)^m (ct)^{2m}}{(2m)!(i)^{2k-2m}} i\sinh{(ct)} \\
 &+ \sum^{k-1}_{m=0}\frac{(-1)^m (ct)^{2m+1}}{(2m+1)!(i)^{2k-2m-1}} \cosh{(ct)} \Bigg], \ \ \ \text{when} \ \  |x| \ even 
\end{cases}$$}

Finally, we can compute the inverse Fourier transform of $\hat{u}$ as \\
\begin{align*}
u(t,x) & = \mathcal{F}^{-1} \Bigg\{ \hat{f}(z)\cos{(czt)} + \hat{g}(z)\frac{\sin{(czt)}}{cz} \Bigg\} \\
&= \mathcal{F}^{-1} \lbrace \hat{f}(z)\cos{(czt)} \rbrace \ + \ \mathcal{F}^{-1} \Bigg\{\hat{g}(z)\frac{\sin{(czt)}}{cz} \Bigg\} \\
& = \ \frac{1}{2\pi} \oint_{C} \hat{f}(z) \cos{(czt)} \ (1+iz)^x \mathrm{d}z + \frac{1}{2\pi c} \oint_{C} \hat{g}(z) \frac{\sin{(czt)}}{z} \ (1+iz)^x \mathrm{d}z \\
& = f \ast \mathcal{F}^{-1}\lbrace\cos{(czt)}\rbrace(x;0) \ + \ \frac{1}{c} \Bigg( g \ast \mathcal{F}^{-1} \Bigg\{\frac{\sin{(czt)}}{z} \Bigg\}(x;0) \Bigg) \\
& = f \ast \mathcal{F}^{-1}\lbrace G_{c,t}(z) \rbrace(x;0) \ + \ \frac{1}{c} (g \ast \mathcal{F}^{-1} \lbrace H_{c,t}(z) \rbrace (x;0)).
\end{align*}

\section{Solving the wave equation with $x$ in a countable time scale}
When x belongs to a countable time scale $\T$, where $\T = \lbrace x_0,x_1,...\rbrace$ , such that $x_0 < x_1 < ...$, we follow the same procedure as when $x\in\Z$, and arrive at the same results, until we come to the task of calculating the inverse Fourier Transforms of $\cos{(czt)}$ and $\frac{\sin{(czt)}}{cz}$. As $x$ belongs to a more complicated time scale, the solutions prove more difficult to compute.

First, we find the expression for the cosine term. When solving for these expressions, we assume that for $x,x^* \in \T, x < x^*$. Also important to the solution is the following bijection $\pi$ defined as $\pi(x_k) = k$ and $\pi^{-1}(k) = x^k$.
\begin{align*}
    \mathcal{F}^{-1} \lbrace \cos{(czt)} \rbrace (x;x^*) \ &= \ \frac{1}{2\pi} \oint_{C} \cos{(czt)} \ e_{iz}(x;x^*) \mathrm{d} z \ = \ \frac{1}{2\pi} \oint_{C} \frac{\cos{(czt)}}{e_{iz}(x^*;x)} \mathrm{d} z\\
    \ &= \ \frac{1}{2\pi} \oint_{C} \cos{(czt)} \prod_{k=\pi(x)}^{\pi(x^*)-1}\frac{1}{1+\mu(x_k)iz} \mathrm{d} z \ \\
    &= \ \frac{1}{2\pi} \prod_{k=\pi(x_k)}^{\pi(x^*)-1}\frac{1}{i\mu(x_k)}\oint_{C} \cos{(czt)} \prod_{k=\pi(x)}^{\pi(x^*)-1}\frac{1}{z-\frac{i}{\mu(x_k)}} \mathrm{d} z. \
\end{align*}
We can then apply the residue theorem from complex analysis to evaluate the remaining contour integral
\begin{equation*}
    \ = \ i\prod_{k=\pi(x_k)}^{\pi(x^*)-1}\frac{1}{i\mu(x_k)} \sum_{j=\pi(x)}^{\pi(x^*)-1} \mathop{\mathrm{Res}}_{z = \frac{i}{\mu(x_j)}}\frac{\cos{(czt)} \prod_{k=\pi(x) , k \neq j}^{\pi(x^*)-1}\frac{1}{z-\frac{i}{\mu(x_k)}}}{z-\frac{i}{\mu(x_j)}}
\end{equation*}
\begin{equation*}
    \ = \ i\prod_{k=\pi(x_k)}^{\pi(x^*)-1}\frac{1}{i\mu(x_k)} \sum_{j=\pi(x)}^{\pi(x^*)-1} \cos{\left(\frac{cti}{\mu(x_j)}\right)} \prod_{k=\pi(x) , k \neq j}^{\pi(x^*)-1}\frac{1}{\frac{i}{\mu(x_j)}-\frac{i}{\mu(x_k)}}.
\end{equation*}
Next, we find the inverse Fourier Transform of the sine term:
\begin{align*}
    \mathcal{F}^{-1} \Bigg\{ \frac{\sin{czt}}{z} \Bigg\} (x;x^*) \ &= \ \frac{1}{2\pi} \oint_{C} \frac{\sin{(czt)}}{z} \ e_{iz}(x;x^*) \mathrm{d} z \ = \ \frac{1}{2\pi} \oint_{C} \frac{\sin{(czt)}}{z} \frac{1}{e_{iz}(x^*;x)} \mathrm{d} z\\
    \ &= \ \frac{1}{2\pi} \oint_{C} \frac{\sin{(czt)}}{z} \prod_{k=\pi(x)}^{\pi(x^*)-1}\frac{1}{1+\mu(x_k)iz} \mathrm{d} z\\
    \ &= \ \frac{1}{2\pi} \prod_{k=\pi(x_k)}^{\pi(x^*)-1} \frac{1}{i\mu(x_k)}\oint_{C} \frac{\sin{(czt)}}{z} \prod_{k=\pi(x)}^{\pi(x^*)-1}\frac{1}{z-\frac{i}{\mu(x_k)}} \mathrm{d} z. \
\end{align*}
We can then apply the residue theorem from complex analysis to evaluate the remaining contour integral to get
\begin{align*}
     i\prod_{k=\pi(x_k)}^{\pi(x^*)-1}\frac{1}{i\mu(x_k)} &\left[ \sum_{j=\pi(x)}^{\pi(x^*)-1} \mathop{\mathrm{Res}}_{z = \frac{i}{\mu(x_j)}} \Bigg\{ \frac{\sin{(czt)}}{z} \prod_{k=\pi(x) , k \neq j}^{\pi(x^*)-1}\frac{1}{z-\frac{i}{\mu(x_k)}} \Bigg\}
     + \mathop{\mathrm{Res}}_{z = 0} \Bigg\{ \frac{\sin{(czt)}}{z} \prod_{k=\pi(x)}^{\pi(x^*)-1}\frac{1}{z-\frac{i}{\mu(x_k)}}\Bigg\} \right]
\end{align*}

\begin{equation*}
    \ = \ i\prod_{k=\pi(x_k)}^{\pi(x^*)-1}\frac{1}{i\mu(x_k)} \sum_{j=\pi(x)}^{\pi(x^*)-1} (-i\mu(x_j))\sin{\Big(\frac{cti}{\mu(x_j)}\Big)} \prod_{k=\pi(x) , k \neq j}^{\pi(x^*)-1}\frac{1}{\frac{i}{\mu(x_j)} - \frac{i}{\mu(x_k)}}.
\end{equation*}
With these expressions in mind for the inverses in mind, we can apply the inverse Fourier transform to $\hat{u}$ and find \\
\begin{align*}
u(t,x) & = \mathcal{F}^{-1} \Bigg\{ \hat{f}(z)\cos{czt} + \hat{g}(z)\frac{\sin{czt}}{z} \Bigg\} = \mathcal{F}^{-1} \lbrace \hat{f}(z)\cos{czt} \rbrace \ + \ \mathcal{F}^{-1} \Bigg\{\hat{g}(z)\frac{\sin{czt}}{cz} \Bigg\} \\
& = \ \frac{1}{2\pi} \oint_{C} \hat{f}(z) \cos{czt} \ e_{iz}(x;x^*) \dd z + \frac{1}{2\pi c} \oint_{C} \hat{g}(z) \frac{\sin{czt}}{z} \  e_{iz}(x;x^*) \dd z \\ 
& = f \ast \mathcal{F}^{-1}\lbrace\cos{czt}\rbrace(x;x^*) \ + \ \frac{1}{c} \Bigg( g \ast \mathcal{F}^{-1} \Bigg\{\frac{\sin{czt}}{z} \Bigg\}(x;x^*) \Bigg).
\end{align*}

\section{Numerical Examples}
We conclude this paper with two numerical examples of approximate solutions to this equation.

First, we look at a numerical solution to the first treated case, where the spatial coordinate is along the integers.
\begin{example}\label{Ex5.1}
When the space time scale is confined to the integers, we find that:
\begin{equation*}
    u^{\Delta^2_{xx}} = u(t,x+2) - 2u(t,x+1) + u(t,x)
\end{equation*}
We can also approximate the second order time derivative with a second order central difference as follows:
\begin{equation*}
    \pd{u}{t} = \frac{u(t+h,x) - 2u(t,x) + u(t-h,x)}{h^2}
\end{equation*}
We can then substitute these expressions into the original wave equation, and solve for the next u value: $u(t+h,x)$.
\begin{equation*}
    u(t+h,x)= c^2 h^2 [u(t,x+2) -2u(t,x+1)] + (2 + c^2 h^2) u(t,x) - u(t-h,x)
\end{equation*}

This equation works for all subsequent time values, except for the first step where we lack information about the point $u(t-h,x)$, which is not defined. We can solve this issue using a trick presented in \cite{numericalbook}. We simply write the central difference approximation for the derivative at t=0, while keeping in mind that $u_t(0,x) = g(x)$. We then solve for $u(t-h,x)$ and substitute the expression into our equation to obtain the following expression for $u(h,x)$.
\begin{equation*}
    u(h,x) = \frac{c^2 h^2}{2} [u(0,x+2) - 2 u(0,x+1)] + \left(1 + \frac{c^2 h^2}{2}\right) u(0,x) + g(x)h
\end{equation*}
Note that $u(0,x) = f(x)$.

In Figure 1 above, we can see an example that was considered for the integers between 0 and 50. At values of x outside of this range, the function was considered to be 0. Without setting an end point for our interval - past which the function takes a value of 0 - the recursion would continue forever, as values of the function at higher x values are always needed.

\begin{figure}
    \centering
    \includegraphics[width=\textwidth]{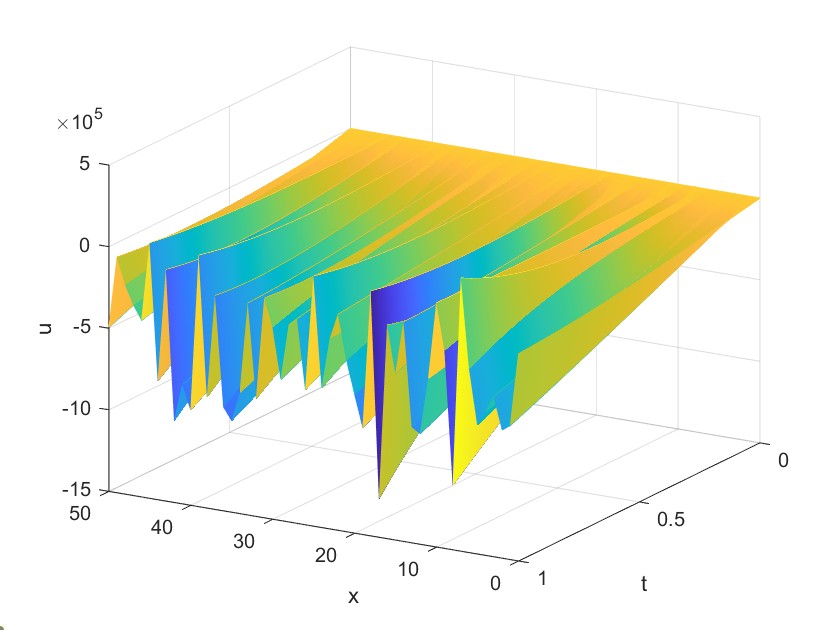}
    \caption{Example of wave equation with $\T = \Z$ and a time step, $h = 0.001$. Lines are propagated in the $t$ direction for each $x$ value, and curves are linearly interpolated between these points. The initial displacement function, $f(x)$, is randomly sampled between $0$ and $1000$ based on a uniform distribution, and the initial velocity function, $g(x)$, is $0$.}
    \label{fig:Zcase}
\end{figure}
\end{example}

For our second example, we look at the case where the spatial coordinate is restricted to the quantum numbers for a positive exponent, which is given as follows.
\begin{equation*}
    q^{\N_0} = \lbrace q^{n} | n \in \N_0 \rbrace \ , \ q>1
\end{equation*}

\begin{example}
When the spatial coordinate is confined to this time scale, we use the fact that on the quantum numbers $\sigma(t) = qt$ and $\mu(t) = t(q-1)$ to find that:
\begin{equation*}
    u^{\Delta^2_{xx}} = \frac{u(t,q^2 x) - (1+q)u(t,qx) + qu(t,x)}{q x^2 (q-1)^2}
\end{equation*}
We once again approximate the second order time derivative with a second order central difference as in Example \ref{Ex5.1}. We can then substitute these expressions into \ref{Eq1.2}, and solve for $u(t+h,x)$ to get
\begin{equation*}
    u(t+h,x) = \frac{c^2 h^2}{qx^2(q-1)^2} [u(t,q^2 x) \ - \ (1+q) u(t,qx)] + \left(2 + \frac{c^2 h^2}{x^2(q-1)^2}\right) u(t,x) - u(t-h,x)
\end{equation*}

Like with the previous example, we have to use an approximation for the first partial derivative with respect to time in order to find the expression for $t = 0$.

\begin{equation*}
    u(h,x) = \frac{c^2 h^2}{2qx^2(q-1)^2} [u(0,q^2x) \ - \ (1+q) u(0,qx)] + \left(1 + \frac{c^2 h^2}{2x^2(q-1)^2}\right) u(0,x) + g(x)h
\end{equation*}
Recall that $u(0,x) = f(x)$.

Figure 2 below shows a graph of the numerical simulation for the quantum numbers. 

\begin{figure}[H]
    \centering
    \includegraphics[width=\textwidth]{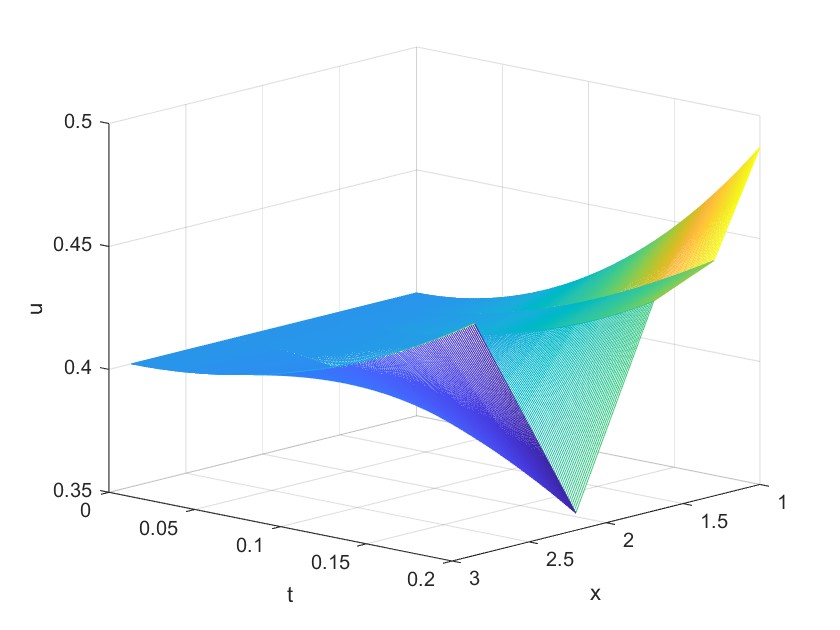}
    \caption{Example of wave equation with $\T = q^{\N_0}$ and a time step, $h = 0.001$. Note that lines are propagated in the $t$ direction for each $x$ value, and curves are linearly interpolated between these points. The initial displacement function, $f(x)$, is a constant, $0.4$, and the initial velocity function, $g(x)$, is $0$. The case with $q = 1.3$ was considered.}
    \label{fig:Qcase}
\end{figure}
\end{example}

\section{Conclusions}
In this paper, solutions were found to the wave equation when the spatial coordinate was confined to the integers and a broad class of discrete time scales. These solutions could be taken a step further by considering the case of a discrete time scale with a non-injective graininess, which would involve calculations with higher order poles. A natural next step would be to examine other well understood partial differential equations, such as the Laplace equation, and try and solve them by similar methods.

Future work includes the case where both the temporal and spatial dimensions of the problem are confined to discrete time scales. Another direction of interest is multi-dimensional formulations of the wave equation in either rectangular and/or polar coordinates. Solutions in higher dimensions will likely necessitate the use of discrete and time scales formulations of special functions and their equations that normally arise from solving partial differential equations in non-rectangular coordinates. Some special functions, like the Bessel function, are better understood in their discrete form \cite{besselpaper}, while others may require work to derive novel formulations.
\vspace{8mm}\\
\textbf{Acknowledgements:}
The authors would like to thank Dr. Tom Cuchta for his help with editing and revising this paper, the insight he provided during the writing process, and for the continued guidance he has provided as a mentor of the first author over the past years.

\bibliography{WaveEqn}
\bibliographystyle{abbrv}

\end{document}